\documentclass[english]{amsart}
\usepackage{babel,enumerate,graphicx}
\usepackage[latin1]{inputenc}

\newtheorem{teor}{Theorem}
\newtheorem{lema}{Lemma}
\newtheorem{prop}{Proposition}
\newtheorem{coro}{Corollary}
\theoremstyle{definition}

\newtheorem{conj}{Conjecture}

\theoremstyle{remark}
\newtheorem{obse}{Remark}

\newtheorem{afir}{Claim}

\newcounter{Casos}

\newcommand{\field}[1]{\mathbb{#1}}

\newcommand{\R}{\field{R}}

\newcommand{\Z}{\field{Z}}

\newcommand{\Cr}[1][1]{\mathrm{C}^{#1}}
\newcommand{\Diff}[1][1]{\mathrm{Diff}^{#1}(M)}
\newcommand{\La}[1][n]{\Delta_{#1}}
\newcommand{\loc}{\mathrm{loc}}
\newcommand{\borrar}[1]{}

\DeclareMathSymbol{\flecha}{\mathbin}{AMSa}{"4C}

\DeclareMathOperator{\dist}{dist}

\begin{document}

\title{A complement to the connecting lemma of Hayashi}

\author[J.C. Mart\'{\i}n]{Jos\'e Carlos Mart\'{\i}n{\dag}}
\address{Departamento de Matem\'aticas Puras y Aplicadas\\
Universidad Sim\'on Bol\'{\i}var\\
Apartado Postal 89000\\
Caracas 1086--A\\
Venezuela}
\email{jmartin@usb.ve}
\thanks{\dag Work partially supported by IMPA-Brasil}
\author[L. Mora]{Leonardo Mora{\ddag}}
\address{Departamento de Matem\'atica\\
Facultad de Ciencias, La Hechicera\\
ULA, Mérida 5101\\
Venezuela}
\email{lmora@ula.ve}
\thanks{\ddag Work partially supported by IMPA-Brasil,
CONICIT-Venezuela and Fundaci\'on Polar.}
\dedicatory{To J. Palis on his 60th anniversary.}
\date{\today}

\begin{abstract} 
Here we show that for a $\Cr[2]$ surface diffeomorphism that satisfy the
hypothesis
of Hayashi connecting lemma either can be approximated, in the 
$\Cr$ topology, by a diffeomorphism exhibiting a homoclinic tangency
or the diffeomorphism already presented transversal homoclinic orbits.
\end{abstract}

\maketitle

\section{Introduction}

\label{sec:intro}
Let $M$ be a $\Cr[\infty]$ compact surface  and  $f:M\rightarrow M $ be a $\Cr$
diffeomorphism. Let $\Lambda$ be an isolated hyperbolic set, we shall
say that $q$ is a {\em homoclinic\/} point (of  $f$)  associated to $\Lambda$
if $q\in
W^s(\Lambda,f)\cap W^u(\Lambda,f)\setminus \Lambda$, where
\begin{gather*}
W^{s}(X,f) = \{x\in M: \lim_{n\to\infty} \dist(f^n(x),f^n(X))=0\}\\
\bigl(\text{resp. $W^{u}(X,f) = \{x\in M: \lim_{n\to -\infty}
\dist(f^n(x),f^n(X))=0\}$}\bigr)
\end{gather*}
is the {\em stable (resp. unstable) set (manifold)} of $X$. When $f$
is well understood from context we shall put simply  $W^{s}(X)$ and
  $W^u(X)$. Let  $p$ be a hyperbolic periodic point of $f$ and $q$ a
homoclinic point associated to $p$, we shall say that $q$ is a {\em
  transversal} homoclinic point if  $T_q W^s(p) + T_q W^u(p)=T_q M$, in
the opposite case we shall say that $q$ is a homoclinic {\em tangency}
associated to $p$. A sequence $\{\gamma_n\}$ of finite orbits
of $f$ is called an    {\em almost homoclinic sequence }
associated to $\Lambda$ if there exist $U$, a neighborhood of $\Lambda$, 
and $\epsilon>0$ such that: 
\begin{enumerate}[i.]
\item $\bigcap_{j\in\Z} f^{j}(U) =\Lambda$,
\item there exist $q^s\in W^s_\epsilon (\Lambda)\setminus \Lambda$, $q^u\in
W^u_\epsilon(\Lambda)\setminus \Lambda$, $q_n'\in \gamma_n$ and $m_n'\in
\Z^+$, for all  $n\in\Z^+$, such that  $q^u=\lim_{n\to\infty} q_n'$ and $q^s
= \lim_{n\to\infty} f^{m_n'}(q_n')$,
\item $\{q_n', f(q_n'), \ldots, f^{m_n'}(q_n')\}\cap (M\setminus U) \neq
\emptyset$, for all $n\in\Z^+$,
\end{enumerate}
where  $W^s_\epsilon(\Lambda)=\bigcup_{x\in\Lambda} W^s_\epsilon(x)$ and
$W^u_\epsilon(\Lambda) = \bigcup_{x\in\Lambda} W^u_\epsilon(x)$ are
the local stable and unstable manifolds of $\Lambda$ respectively.
We  shall say  that an invariant  set $\Lambda$  of $f$ is a 
{\em basic set } if it is a compact, hyperbolic, isolated and transitive.

  The $\Cr[1]$ connecting lemma of Hayashi  \cite{h} says that, if a
compact isolated hyperbolic set  $\Lambda$ has associated  an almost homoclinic
sequence, then for all  $\Cr[1]$ neighborhood $\mathcal{U}$ of $f$,
there exists  $g\in \mathcal{U}$  coinciding with $f$ in a neighborhood 
of $\Lambda$  and having a homoclinic point associated to a periodic
hyperbolic
point in $ \Lambda$.
 Our purpose here is  to show the following theorem. 

    \begin{teor}    
If  $f:M^2\rightarrow M^2$ is a $\Cr[ 2]$ diffeomorphism with a   basic set
 $\Lambda$ which has associated an almost homoclinic sequence
then for every periodic point $p \in \Lambda$ one of the following statements hold: 
\begin{enumerate}[(i)] 
\item $p$ has associated a transversal homoclinic point outside $\Lambda$,
\item for all  $\mathcal{N}\subset \Diff[ 1]$, a neighborhood of $f$,
 there exists
  $g\in\mathcal{N}$
having a homoclinic tangency associated to $p$.
\end{enumerate}
\end{teor}

In \cite{ps}, it is proved that given a diffeomorphism
$f:M^2\rightarrow M^2$, it can be $\Cr[1]$
approximated by a diffeomorphism exhibiting a homoclinic tangency or
by an Axiom A diffeomorphism. Since the dynamical richness exhibited
by nearby diffeomorphisms of  another one exhibiting a homoclinic
tangency \cite{pt}, 
 it is an important problem   to decide when,  in a particular
 dynamical situation, the system can be approximated by another one
 having a homoclinic
 tangency. We raise this question in the following dynamical
 situation: Suppose that $f:M^2\rightarrow M^2$ is a $\Cr[ 2]$ diffeomorphism
 which has a compact basic  set $\Lambda$ and $q \in ((\overline{W^s (\Lambda)}
 \cap W^u (\Lambda)) \cup (\overline{W^u (\Lambda)}
 \cap W^s (\Lambda))) \setminus \Lambda$. The Hayashi $\Cr[1]$ connecting lemma
 says to us that in this case we can approximate $f$ by another
 diffeomorphism
$g$ having a homoclinic orbit. But it does not say what kind of
homoclinic orbit it has. The following corollary of theorem
1  give us the answer.

\begin{coro} Let $f:M^2\rightarrow M^2$ is a $\Cr[ 2]$ diffeomorphism
 which has a  basic set  $\Lambda$ and $q \in ((\overline{W^s (\Lambda)}
 \cap W^u (\Lambda)) \cup (\overline{W^u (\Lambda)}
 \cap W^s (\Lambda))) \setminus \Lambda$.  If $W^s (\Lambda) \cap W^u
 (\Lambda)\setminus\{\Lambda\}=\emptyset$, then for every $\Cr[1]$ neighborhood
$\mathcal{N}\subset \Diff[ 1]$ of
$f$ and for every  hyperbolic periodic point $p \in \Lambda$,
there exists $g \in \mathcal{N}$ having a homoclinic tangency
 associated to a $p$.
 \end{coro}
We do not know if the $\Cr[2]$ hypothesis on $f$ above can be
raised. The reason of its appearance is that we use  results of \cite{ps}
which assume this hypothesis.

Respect to what happens in higher dimensions, we remark that examples 
of Díaz \cite{d} make the theorem 1 false. What we think 
should be true is the following.

\begin{conj}  
 If  $f:M \rightarrow M $ is a $\Cr[ 2]$ diffeomorphism with a compact basic
 set $\Lambda$ which has associated an almost homoclinic sequence,
 then for every hyperbolic
periodic point $p\in \Lambda $  one of the following statements hold: 
\begin{enumerate}[(i)] 
\item $p$ has associated a transversal homoclinic point outside  $\Lambda$,
\item for all  $\mathcal{N}\subset \Diff[ 1]$, a neighborhood of $f$,
 there exists
$g\in\mathcal{N}$
having a homoclinic bifurcation associated to $p$.
\end{enumerate} 
\end{conj}
Here we use the notion of homoclinic bifurcation as defined in
\cite{pt}, page 134.

The paper is organized in three sections. In the second one we 
present all stuff which is necessary for doing the proof of the 
theorem, and in the last section we provide the proof of the theorem.
The proof's idea  is the following. Let $f$ and $\Lambda$ be as in the
theorem. By  the  connecting lemma of
Hayashi, $f$ can be approximated by diffeomorphisms exhibiting a
transversal homoclinic orbit.  Analyzing the rate 
of  domination between stable and unstable vectors for these homoclinic
orbits, we obtain:  If these rates are bounded
for infinitely many  perturbations, then $f$ has  a homoclinic orbit
outside $\Lambda$.  In  the other case $f$ can be 
approximated by diffeomorphisms exhibiting a homoclinic tangency.
   
\section{Preliminaries}

Consider  a diffeomorphism $f: M\rightarrow M $ and a  basic set
$\Lambda$ of $f$ which has associated an almost homoclinic sequence $\{\gamma_n\}$.
Let $q^s$, $q^u$, $q_n'$ and  $m_n'$ as in definition of an almost
homoclinic sequence. Denote by $p^s$ and $p^u$ the points of $\Lambda
$ such that $q^s\in W^s_{\loc}(p^s)$
and  $q^u\in W^u_{\loc}(p^u)$. Let  $p\in \Lambda$ be any  hyperbolic
periodic point.

%%%%%%%%%%%%%%%%%%%
\begin{figure}
\begin{center}
\includegraphics*{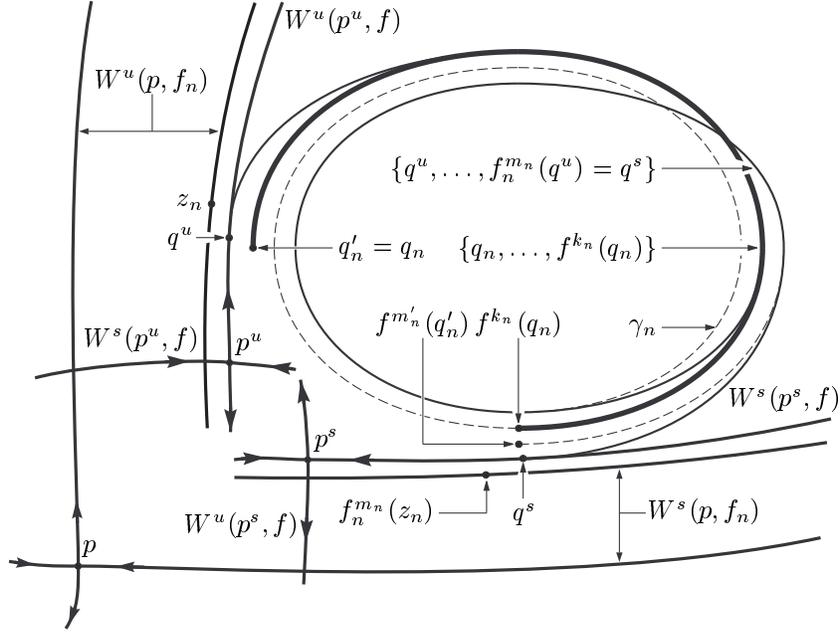}
\caption{Proof of Lemma \ref{lem:conlem}: Connection between stable and
  unstable manifolds of the orbits of  $p^s$ and  $p^u$, respectively.
Compare with  Figure~1(B) of  \cite{h}.}
\label{fig1}
\end{center}
\end{figure}
%%%%%%%%%%%%%%%%%%%

\begin{lema}
\label{lem:conlem}
If $W^s(\Lambda)\cap W^u(\Lambda)\setminus \Lambda =\emptyset$ then 
there exist sequences  $\{f_n\}\subset \Diff[ 1]$, $\{z_n\}\subset M$,
$\{q_n\}\subset M$, $\{k_n\}\subset \Z^+$ and  $\{m_n\}\subset \Z^+$
such that:
\begin{enumerate}[i.]
\item\label{lemit:1} $f_n\to f$, in the  $\Cr[ 1]$topology,
\item for each $n$, $z_n$ is a transversal homoclinic point of $f_n$  associated  to $p$,

\item\label{lemit:3} $\lim_{n\to\infty} z_n=q^u$ and $\lim_{n\to\infty}
f_n^{m_n}(z_n) = q^s$;
\item\label{lemit:4} $\lim_{n\to\infty} q_n=q^u$ and $\lim_{n\to\infty}
f^{k_n}(q_n) = q^s$, and
\item\label{lemit:5} $\forall\, \epsilon>0$ there exists $N\in\Z^+$
  such that 
$f^j(q_n)\in\bigcup_{k=0}^{m_n} B_\epsilon(f_n^k(z_n))$, for all $0\le
j\le k_n$ and  $n\ge N$.
\end{enumerate}
\end{lema}

\begin{obse}
\label{obs:1}

We observe that the hypothesis $W^s(\Lambda)\cap W^u(\Lambda)\setminus \Lambda=
\emptyset$ is equivalent to $m_n'\to \infty$ for all almost homoclinic
sequences.  Because if $m = \liminf_{n\to\infty} m_n'<\infty$ then 
$f^m(q^u)=q^s$ that implies that $W^s(\Lambda)\cap W^u(\Lambda)\setminus
\Lambda\neq \emptyset$. The converse is obvious.
\end{obse}

For the proof of this lemma we use the following version of the
connecting lemma. 
\begin{teor}
\label{teo:wx}
Let  $M$ be a  Riemannian manifold, $f:M\rightarrow M$ be a $\Cr$
 diffeomorphism
 and  $z$
be a non periodic point of $f$. Then for every neighborhood
 $\mathcal{U}$ of $f$ in 
$\Diff$ there exist  $\rho>1$, $L\in\Z^+$ and $\delta_0>0$ such that 
for each 
$0<\delta\le \delta_0$ and each couple of points, $x\notin T= \bigcup_{n=1}^L
f^{-n} B_\delta(z)$ and $y\in B_{\delta/\rho}(z)$, for which the  future orbit 
of  $x$ mets  $B_{\delta/\rho}(z)$, there exist 
$g\in\mathcal{U}$ and  $m\in \Z^+$ such that  $g=f$ outside of  $T$
 and  $g^m(x)=y$.
\end{teor}
This theorem is  Theorem F of 
\cite{wx}.

\begin{obse}
\label{obs:2} In the conditions and notation of 
theorem~\ref{teo:wx}. If $k\ge L$ is the first positive integer such 
that  $f^k(x)\in B_\delta(z)$, then  $k\le m$,
$f^k(x)\in B_{2\delta}(y)$ and $f^j(x)\in \bigcup_{n=0}^m
B_{2\delta}(g^n(x))$, for all $0\le j\le k$.
\end{obse}
\begin{proof}[Proof of lemma 1]

Because of the remark \ref{obs:1}, we know that $m_n'\to \infty$.
Take a sequence  $\{\epsilon_n\}$ such that  $\epsilon_n\to 0$. For
each $n$
we take  $\rho_n^-$, $L_n^-$ and  $\delta_n^-$ as in  Theorem~\ref{teo:wx} 
for $\mathcal{U}=B_{\epsilon_n}(f) = \{g\in\Diff: \dist_{\Cr}(g,f) <
\epsilon_n\}$ and constants $\rho_n^+$, $L_n^+$ and $\delta_n^+$ for
$B_{\epsilon_n}(f^{-1})$. Let $\rho_n = \max \{\rho_n^-, \rho_n^+\}$, $L_n
= \max \{L_n^-, L_n^+\}$ and  $\delta_n = \min \{\delta_n^-, \delta_n^+\}$.
Given  $\delta>0$ we define  $T_n^-(\delta) = \bigcup_{j=1}^{L_n}
f^{-j}(B_{\delta}(q^s))$ and  $T_n^+(\delta) = \bigcup_{j=1}^{L_n}
f^{j}(B_{\delta}(q^u))$. It is easy to see there exist $0<\delta_n'\le \delta_n$
such that
\begin{enumerate}[i.]
\item $\{B_{\delta_n'}(q^u), f(B_{\delta_n'}(q^u)),\ldots ,
f^{2L_n}(B_{\delta_n'}(q^u)), B_{\delta_n'}(q^s), \ldots,
f^{-2L_n}(B_{\delta_n'}(q^s))\}$ is a  family of disjoint  sets,
\item $\overline{T_n^{\pm}(\delta_n')} \cap \Lambda = \emptyset$,
\item $T_n^-(\delta_n') \cap \bigcup_{j\ge 0} f^{-j}(W^u_{\loc}(p^u)) =
\emptyset$ and $T_n^+(\delta_n') \cap \bigcup_{j\ge 0} f^{j}(W^s_{\loc}(p^s))
= \emptyset$.
\end{enumerate}
In particular, $T_n^-(\delta_n')\cap T_n^+(\delta_n')= \emptyset$. Let
$\{\alpha_n\}\subset \Z^-$ such that  $q_{\alpha_n}'\in
B_{\delta_n'/\rho_n}(q^u)$ and  $f^{m_{\alpha_n}'}(q_{\alpha_n}')\in
B_{\delta_n'/\rho_n}(q^s)$, for easy of writing let us assume that
$\alpha_n=n$, see figure ~\ref{fig1}. By theorem~\ref{teo:wx} there
are a sequence of diffeomorphisms  $\{g_n\}$ and a sequence of
positive integers 
$\{m_n\}$ such that  $g_n = f$ outside of  $T_n^-(\delta_n')$, $g_n^{m_n}(q_{
n}') = g_n^{m_n-L_n-1} (f^{L_n+1} (q_{ n}' ) ) = q^s$. In the same way,
we get a sequence of diffeomorphisms  $\{\tilde{g}_n\}$ such that 
$\tilde{g}_n = f$ outside of  $T_n^+(\delta_n')$ and  $\tilde{g}_n^{L_n+1}(q^u) =
f^{L_n+1}(q_{ n}')$.

Let  $f_n:M\rightarrow M $ be defined by 
$$
f_n(x) = \begin{cases} f(x) & \text{if $x\notin T_n^-(\delta_n')\cup
T_n^+(\delta_n')$}\\
g_n(x) & \text{if  $x\in T_n^-(\delta_n')$}\\
\tilde{g}_n(x) & \text{if $x\in T_n^+(\delta_n')$.}\end{cases}
$$
Because of the construction,  $f_n \to f$,  $\Lambda$ is a hyperbolic
set of  $f_n$, $f_n(q^u) = q^s$, $W^s_{\loc}(p^s,f_n) =
W^s_{\loc}(p^s,f)$ and  $W^u_{\loc}(p^u,f_n) = W^u_{\loc}(p^u,f)$.
We can suppose that  $W^s(p^s,f_n)$ and  $W^u(f^{m_n}(p^u), f_n)$ meet
in a transversal way at  $q^s$, if not with a small perturbation we get
that.
 Since  $\bigcup_{j=0}^{k-1} W^s(f^j(p),f_n)$ and 
$\bigcup_{j=0}^{k-1} W^u(f^j(p), f_n)$ are dense in  $\Lambda$, for all
$n$, where  $k\ge 1$ is the period of   $p$, we have that $q^u$ 
is accumulated by transversal homoclinic points of  $f_n$ associated to
$p$, so   $q^u$ is accumulated by transversal homoclinic points of  $f_n$ associated to
$p$. In particular we can take  $z_n$, a homoclinic point of $f_n$
associated to  $p$, in such a way that it holds 
\eqref{lemit:3}. The existence of  $q_n$ and  $k_n$ satisfying 
\eqref{lemit:4} and  \eqref{lemit:5} is an immediate consequence of
the previous thing and of remark~\ref{obs:2}.
\end{proof}

The next lemma is the basic tool for making perturbations and is a 
slightly general version of \cite[Lemma 1.1]{f}.

\begin{lema}
\label{lem:f}
Let  $f\in \Diff[ 1]$. Given  $\mathcal{U}$, a $\Cr[ 1]$
neighborhood of $f$ , there exist  $\epsilon >0$ and a neighborhood
$\mathcal{U}_0\subset \mathcal{U}$ of $f$ such that the following
holds.
 If we have   $\tilde{f}\in\mathcal{U}_0$, $\{x_1,\ldots, x_n\}\subset M$, $U$
 a neighborhood of  $\{x_1,\ldots, x_n\}$ and  $T_i:T_{x_i}M\to
T_{\tilde{f}(x_i)}M$  linear mappings with
$\|T_i-D_{x_i}\tilde{f}\|< \epsilon$, for all $1\le i\le n$, then
there exists $g\in\mathcal{U}$ such that $g(x)=\tilde{f}(x)$, for
$x\in \{x_1,\ldots, x_n\}\cup (M\setminus U)$, and  $D_{x_i}g = T_i$, for
all  $1\le i\le n$.
\end{lema}

The following two lemmas are the tools to perturb $f$ along a piece of
a homoclinic orbit to close the angle made by the stable and unstable
manifolds.

\begin{lema}
\label{lem:1}
For each $\delta>0$ and $\theta>0$ there exists $c>1$ such that the
following holds. If $|b_1|/|b_2|\ge c$ and
$$
B=\begin{bmatrix} b_1 & 0\\ 0& b_2\end{bmatrix},
$$
then the angle between  $e_1$ and $B\circ I_\delta e_2$ is less than
$\theta$, where $I_\delta = \bigl[\begin{smallmatrix} 1 & \delta\\ 0 &
  1\end{smallmatrix}\bigr]$, $e_1=
(\begin{smallmatrix}
  1 \\ 0
\end{smallmatrix})
$ and  $e_2 =( \begin{smallmatrix}
  0 \\ 1
\end{smallmatrix})$.
\end{lema}

\begin{proof}
We may assume that 
 $\theta < \frac{\pi}{2}$.  If  $c>\frac{1}{\delta\tan \theta}$, then $\angle (e_1, B\circ I_\delta e_2) \le \angle (e_1, e_1+c\delta e_2) < \theta$.
\end{proof}

Now consider a sequence $\{A_n\}$ of diagonal matrices

$$
A_n = \begin{bmatrix} a_n & 0\\ 0 & b_n\end{bmatrix},
$$
with  $a_n, b_n>0$. Let  $\sigma_n = \prod_{j=0}^{n-1} a_n$ and $\tau_n = \prod_{j=0}^{n-1} b_n$.

\begin{lema}
\label{lem:2}
For each  $\epsilon > 0$, $\theta>0$ and  $K>0$ there exists
$r_0\in\Z^+$ such that the following holds. If $\{A_n\}$ is a sequence
of matrices as above with $\|A_n\|\le K$ and $\sigma_r/\tau_r > 1/2$,
for $r\ge r_0$, then there exist  matrices $\tilde{A}_n$, for
$n=0,\ldots ,r-1$, such that  $\|\tilde{A}_n - A_n\|<\epsilon$,
$\tilde{A}_n e_1 = e_1$ and the angle between  $e_1$ and
$\tilde{A}_{r-1}\circ \cdots \circ \tilde{A}_0 e_2$ is less than  $\theta$.
\end{lema}

\begin{proof}
First of all, we observe that there exists $\delta = \delta(\epsilon,
 K)$
such that if  $C$ is a matrix $\delta$-close to the identity and $D$
 is a matrix with  $\|D\|\le K$ then $\|C\circ D-D\|<\epsilon/2$.

Let  $c>1$ be given by  Lemma~\ref{lem:1} for $\delta$ and $\theta$,
and let $r_0\in\Z^+$ such that  $(1+\delta)^{2r_0}\ge 2c$. If $\{A_n\}$ is a
sequence as in the hypothesis and $r\ge r_0$, we define 
$$
\hat{A}_n = \begin{bmatrix} (1+\delta)a_n & 0\\ 0 &
  \displaystyle{\frac{1}{1+\delta}} b_n\end{bmatrix},\quad \text{for  $n=0,\ldots ,r-1$}.
$$
Clearly $\hat{A}_n$ is $\epsilon/2$-close to $A_n$. In addition, $B=
\hat{A}_{r-1}\circ \cdots \circ \hat{A}_0$ satisfies the hypothesis of
 Lemma~\ref{lem:1}, so the angle between  $e_1$ and $B\circ I_\delta e_2$
is less than   $\theta$. Put $\tilde{A}_n = \hat{A}_n$, for $1\le n\le
r-1$, and $\tilde{A}_0 = \hat{A}_0\circ I_\delta$.
These matrices fulfill the thesis, so we are done. 
\end{proof}

Now we introduce some notation before presenting the next lemma.
Let $x\in M$ be a fixed point of $f$, $\sigma$ and $\tau$ be the
eigenvalues of  $D_xf$ with  $|\sigma|<|\tau|$. By the center stable 
manifold theorem, there exist invariant manifolds $W^-(x)$ and $W^+(x)$
associated to $\sigma$ and  $\tau$ respectively. We shall say
that $\Sigma$ is a {\em separatrix\/} of  $W^\pm(x)$ if $\Sigma$
is one of the connected component of $X\setminus \{x\}$, where
$X$ is the connected component of $W^\pm(x)\cap (W^s(x)\cup W^u(x))$
which contains $x$. A separatrix $\Sigma$ will be called 
{\em stable\/} (resp. {\em unstable\/}) if  $\Sigma\subset W^s(x)$
(resp. $\Sigma\subset W^u(x)$). When $x$ is a periodic point
same definitions for separatrices extend considering  $f^k$, where
$k$ is the period of  $x$. Observe that it can happen that a periodic 
point $x$ has not separatrices of $W^\pm(x)$ at all. Denote by
$\ell(\Sigma)$ the length of  $\Sigma$.

Let $\Delta\subset M$ be an invariant compact set of $f$ and 
$E\oplus F=T_{\Delta}M$ be a  $Df$-invariant splitting. We shall
say that   $E\oplus F$ is a {\em dominated splitting\/}  of $\Delta$
if there exist $C>0$ and  $0<\lambda<1$ such that
$$
\left\|Df^j|_{E(x)}\right\| \left\|Df^{-j}|_{F(f^j(x))}\right\| \le C\lambda^j,
$$
for all  $j\in\Z^+$ and  $x\in\Delta$.

\begin{lema}
\label{lem:3}
Let  $\Delta$ be a compact invariant set with a
 dominated splitting $E\oplus F$, $x\in \Delta$ be a periodic point
, $\Sigma$ be an unstable (resp. stable) separatrix  of  $W^+(x)$
 (resp. $W^-(x)$) and $y\in \Sigma\cap \Delta$. If
 $\ell(\Sigma)<\infty$ then there exist a neighborhood  $U$ of $y$ and
a periodic point  $z$ such that  $U\subset \bigcup_{n=1}^k W^s(f^n(z))$
 (resp. $U\subset \bigcup_{n=1}^k W^u(f^n(z))$), where  $k$ is the
 period of  $z$.
\end{lema}

\begin{proof} We assume that   $x$ is a fixed point, if $x$ is a
 periodic
point take $f^k$, where $k$ is the period of $x$.
Let us  assume  that  the eigenvalues of $D_xf$ are positive.
We prove the result when $\Sigma$ is an unstable separatrix, the other
cases
being analogous. As $\ell(\Sigma)<\infty$, the limit
$\lim_{n\to\infty} f^n(y)$ exists. Let $z$ denote this limit, then  
$z \in \Delta$. Let  $\sigma$ and  $\tau$ be the eigenvalues  of
$D_zf$, with $|\sigma|\le |\tau|$. Note that  $|\sigma| < |\tau|$
since $z \in \Delta$ and $|\sigma|\le 1$.

%%%%%%%%%%%%%%%%%
\begin{figure}
\begin{center}
\includegraphics*{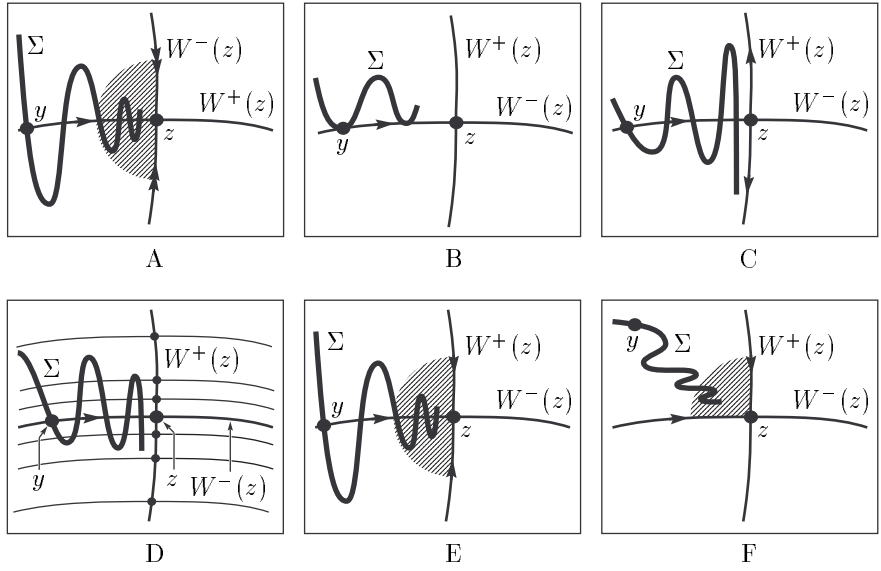}
\caption{Proof of  Lemma~\ref{lem:3}. 
The shadow parts correspond to a piece of the stable manifold of  $z$.}
\label{fig2}
\end{center}
\end{figure}
%%%%%%%%%%%%%%%%

To start, assume that $y$ belongs to any of invariant manifolds
of $z$. If $y\in W^+(z)$, see Figure~\ref{fig2}-A, then $|\sigma|<1$
and consider two cases: (1) $|\tau|<1$ and  (2) $|\tau|=1$. In the
first case  $z$ is a sink and the result follows. In the second case 
there exists a strong contractive foliation in a neighborhood of  $z$,
so the dynamics near to $z$ is decided by the dynamics on  $W^+(z)$.
As $y$ belongs to a stable separatrix of $z$ we are done.

Now assume that $y\in W^-(z)$. Observe that at $y$, $\Sigma$ and
$W^-(z)$
meet transversally, because if the intersection were non transversal
then $T_xW^+(x)=F(x)$, $T_zW^-(z)=E(z)$ and  $y\in W^+(x)\cap
W^-(z)\cap \Delta$, see  Figure~\ref{fig2}-B, would imply that
$E(y)=F(y)$, which  is an absurd. Consider the following three cases:
(1) $W^+(z)$ has an unstable separatrix, (2) $W^+(z)$ has no two
separatrices and (3) both separatrices of $W^+(z)$ are stable.
In the first case we have $\ell(\Sigma)=\infty$, which follows from
the $\lambda$-lemma when  $|\tau|>1$, and from the existence  of a 
strong contractive foliation when  $|\tau|=1$, see  Figure~\ref{fig2}-C. 
In the second case we get also the same conclusion,  observe that here
$|\tau|=1$ and $z$ is accumulated by fixed points in   $W^+(z)$, see
Figure~\ref{fig2}-D. In the third case the result follows
easily as seen in  Figure~\ref{fig2}-E.

If $y$ does not  belong to $W^-(z)\cup W^+(z)$, then using analogous 
arguments as above we get the existence of $U$, see
Figure~\ref{fig2}-F.

If  any of the eigenvalues of $D_xf$ is negative, the proof runs as
above, but in this case it is possible that  $z$ could be a period point 
of period two.
\end{proof}

\section{Theorem's Proof}

Let   $\Lambda $ be a  basic set, $\{\gamma_n\}$ be an almost
homoclinic sequence associated to $\Lambda$ and  $p\in \Lambda$ be a
periodic point. To simplify, we suppose that $p$ is a fixed point.
 Let  $U$, $q^s$,
$q^u$, $\{q_n'\}$ and  $\{m_n'\}$ as in the definition of an almost
homoclinic sequence.  We can assume that  $\{\gamma_n\} = \{q_n',\ldots
,f^{m_n'}(q_n')\}$.

Let us see that we get the result when $W^s(\Lambda)\cap
W^u(\Lambda)\setminus \Lambda\neq \emptyset$. 
Let  $x^s, x^u\in \Lambda$
such that  $W^s(x^s)\cap W^u(x^u)\neq \emptyset$. 
If there is a transversal meeting point between  $W^s(x^s)$ y
$W^u(x^u)$,
then the result follows from the density of the manifolds
 $W^s(p)$ and  $W^u(p)$ in 
$\Lambda$.
If   $W^s(x^s)$ and  $W^u(x^u)$ does not intersects transversally
then, using similar arguments to the Theorem 1's proof of
\cite[pag. 36]{pt}, we get that  $f$ can be approximated in the 
$\Cr$ topology by a diffeomorphism which has a homoclinic tangency
associated to
$p$.

In what follows, we assume that 
 $W^s(\Lambda)\cap W^u(\Lambda)\setminus
\Lambda = \emptyset$.
Let  $\{f_n\}$, $\{z_n\}$, $\{q_n\}$, $\{m_n\}$ and
$\{k_n\}$ as in  lemma~\ref{lem:conlem}, and 
$\Delta_n=\{f_n^k(z_n):k\in\Z\}\cup\{p\}$. It is clear that
 $\Delta_n$ is a hyperbolic set for $f_n$, so it has a 
dominated splitting 
$E_n\oplus F_n$. By the same argument that in  remark~\ref{obs:1},
$m_n\to \infty$ and $k_n\to\infty$.

Given  $n\in\Z^+$ we denote by  $l_n$ the least positive integer such
that 
$$
\left\|Df_n^k|_{E_n(x)}\right\| \left\|Df_n^{-k}|_{F_n(f_n^k(x))}\right\|
\le \frac{1}{2},
$$
for all  $k\ge l_n$ and  $x\in\Delta_n$.
 Since  $\La$ has a dominated splitting, $l_n$ is well defined.
The proof follows from the consideration of two cases: 1)
$\liminf_{n\to\infty} l_n=\infty$, and  2) $\liminf_{n\to\infty}
l_n=l<\infty$. Taking subsequences, if necessary, we can suppose that 
$l_n\to\infty$ (resp. $l_n = l$) in the first  (resp. second ) case.
The next lemma finish the proof in the first case. 
 
\begin{lema}
If  $l_n\to\infty$ then there exists a sequence  $\{g_n\}\subset
\Diff[ 1]$ such that $g_n\to f$ in the  $\Cr[ 1]$ topology, $p$ is
a hyperbolic fixed point of  $g_n$ and  $g_n$ has a homoclinic
tangency associated to  $p$.
\end{lema}

\begin{proof} It is enough to show that,
given $\theta>0$ and  $\mathcal{U}$ neighborhood of $f$,
there exists $g\in \mathcal{U}$ such that  $p$
is a hyperbolic fixed point of  $g$ and there is a 
homoclinic point of $g$, $x$ associated to $p$, such that 
the angle between $T_xW^s(p)$ and  $T_xW^u(p)$ is less than 
 $\theta$. Let  $\theta>0$ and  $\mathcal{U}$ a neighborhood of
 $f$ in  $\Diff$.

Let $N\in\Z^+$ such that $f_n\in\mathcal{U}$, for all $n\ge N$.
If there exist $x_n\in \Delta_n$ such that $\liminf_{n\to\infty}
\angle(E_n(x_n),F_n(x_n))<\theta$ then we have that we want. So assume 
that $\angle(E_n(x),F_n(x))\ge \theta$ for all  $x\in\Delta_n$ and $n\ge N$.

By hypothesis, there exist $x_n\in \Delta_n$ and  $r_n\in\Z^+$, for
each $n\in\Z^+$, such that $r_n\to\infty$ and 
$$
\left\|Df_n^{r_n}|_{E_n(x_n)}\right\| \left\|Df_n^{-r_n}|_{F_n(f_n^{r_n}(x_n))}\right\| > \frac{1}{2}.
$$
It is clear that $x_n\neq p$ for all  $n$. Let
 $u_n^k\in E_n(f_n^k(x_n))$ and  $v_n^k\in F_n(f_n^k(x_n))$, with
 $\|u_n^k\|=\|v_n^k\|=1$, for all $k,n\in \Z^+$. Let 
$B_n(k)$ be the matrix of $Df_n(f_n^k(x_n))$ in  the base
$\{u_n^k,v_n^k\}$.
Observe that  $B_n(k)$ is a diagonal matrix and there exists 
$K>0$ such that  $\|B_n(k)\|\le K$, for all $k,n\in \Z^+$.
Let  $\epsilon > 0$, given by  Lemma~\ref{lem:f} for  $f$ and
$\mathcal{U}$; $r_0\in\Z^+$, given by Lemma~\ref{lem:2} for
$\epsilon$, $\theta$ and  $K$, and  $n\ge N$ such that $r_n\ge r_0$.
We define $A_k = B_n(k)$, for  $k\ge 0$. Observe that 
$\{A_k\}_{k=0}^{r_n}$ satisfy the hypothesis of Lemma~\ref{lem:2}, so
there are matrices $\tilde{A}_0,\ldots ,\tilde{A}_{r_n-1}$ such that
 $\|\tilde{A}_k-A_k\|<\epsilon$, $\tilde{A}_k e_1 = e_1$
and the angle between  $e_1$ and $\tilde{A}_{r_n-1}\circ \cdots \circ
\tilde{A}_0 e_2$ is less than  $\theta$. 
By Lemma~\ref{lem:f}, there exits  $g\in \mathcal{U}$ such that $g^j(x_n) =
f_n^j(x_n)$, for all $k\in\Z$ and the angle between  $T_{y}W^s(p)$ and 
$T_{y}W^u(p)$ is less than  $\theta$, with $y=f_n^{r_n}(x_n)$, that is
what we wanted to proof. 
\end{proof}

Lets go consider, to finish,  the second case, i.e. when  $l_n=l$. We will
see that this case drives us to a contradiction if we assume that, as
in fact we are doing, that  $W^s(\Lambda)\cap W^u(\Lambda)\setminus \Lambda =
\emptyset$.

There exist a sequence $\{n_i\}$ such that $\{\La[n_i]\}$ converges,
in the Hausdorff topology, to a compact set $\Delta$. 
To simplify  suppose that 
 $n_i=i$. Since $l_n=l$, then  $\Delta$ has a dominated splitting
 $E\oplus F$. The subbundles  $E$ and  $F$ can be extended in a
 continuous way to a neighborhood 
  $V_0$ of $\Delta$ in such a way that  $E(x)\oplus
F(x)=T_xM$, for all  $x\in V_0$. Given  $\eta\in T_x M$, for $x\in V_0$,
we denote by $\eta_1\in E(x)$ and  $\eta_2\in F(x)$ to the vectors 
such that 
$\eta=\eta_1 + \eta_2$. Let 
$$
C^E_r (x)=\{\eta\in T_xM: \|\eta_2\|< r\,\|\eta_1\|\}\quad \text{y}\quad
C^F_r (x)=\{\eta\in T_xM: \|\eta_1\|<r\,\|\eta_2\|\},
$$
for  $r>0$.

Let  $V_1$ be an open neighborhood of  $\Delta$, $0<\rho <1$ and $\frac{1}{2} <
\lambda < 1$ be constants such that  $\overline{V_1}\subset V_0$,
\begin{gather*}
Df^{-1} C^E_\rho(x)\subset C^E_\rho(f^{-1}(x)),\quad Df C^F_\rho(x)\subset
C^F_\rho(f(x))\quad\forall\, x\in \overline{V_1}\\
\bigl\|Df^k|_{C^E_\rho(x)}\bigr\|
\bigl\|Df^{-k}|_{C^F_\rho(f^k(x))}\bigr\|< \lambda \quad\forall\, k\ge l\
\forall\, x\in \bigcap_{j=0}^{k}f^{-j}(\overline{V_1}).
\end{gather*}
Let  $\La[0]^{\pm}=\bigcap_{j=0}^{\infty} f^{\pm j}(\overline{V_1})$ and
$\Delta_0=\La[0]^+\cap \La[0]^-$. We observe that  $\Delta\subset \La[0]$,
$f(\La[0]^+)\subset \La[0]^+$ and $f^{-1}(\La[0]^-)\subset \La[0]^-$.

We call an  {\em interval\/} to the image of an immersion 
in  $M$ of a closed interval and we will denote by  $\ell(I)$ its
length. We shall say that an interval  $I$ is an  {\em $E$-interval\/}
if  $I\subset \overline{V_1}$ and 
$T_xI\subset C^F_\rho(x)$, for all  $x\in I$. Analogously, we define 
{\em $F$-intervals.}  An  $E$-interval $I$ is a 
$\delta$-$E$-{\em interval} if  $I\subset \La[0]^+$ and
$\ell(f^k(I))<\delta$, for all 
 $k\ge 0$. In the same way, we say  that an  $F$-interval $I$ is a {\em
$\delta$-$F$-interval \/} if $I\subset \La[0]^-$ and 
$\ell(f^{-k}(I))<\delta$, for all $k\ge 0$.

Now we state Proposition~3.1 of \cite{ps}.

\begin{prop}
\label{pro:ps}
There exists $\delta_0>0$ such that if  $I$ is a  $\delta$-$E$-interval, with
$0<\delta\le \delta_0$, then one of the following alternatives hold:
\begin{enumerate}[i.]
\item $\omega(I)=\bigcup_{k=1}^m \mathcal{C}_k$, where each  $\mathcal{C}_k$
is a simple closed curve, normally hyperbolic for  $f^m$, with its
dynamics an irrational rotation  and $f(\mathcal{C}_{k-1}) =
\mathcal{C}_k$, for  $1\le
k\le m$, with  $\mathcal{C}_0=\mathcal{C}_m$.
\item $\omega(I)$ has just periodic points.
\end{enumerate}
\end{prop}
We recall that there are points  $q_n$ such that  $q_n\to q^u$ and
$f^{k_n}(q_n)\to q^s$, see Figure~\ref{fig1}. We observe that  $q^s,q^u\in
\Delta$, because of  $q^s = \lim_{n\to \infty} f_n^{m_n}(z_n)$, $q^u = \lim_{n\to
\infty} z_n $ and  $z_n \in \La$, for all $n$. Let  $V_2$ be an open
neighborhood of  $\Delta$, $\delta_1>0$ and  $N_0\in\Z^+$ be such that  $\overline{V_2}\subset
V_1$, $B_{\delta_1}(V_2)\subset V_1$ and  $\La\subset V_2$, for  $n\ge N_0$.
The item  \eqref{lemit:5} of  Lemma~\ref{lem:conlem} implies that we
there exists  $N_1\ge
N_0$ such that  $f^j(q_n)\in V_2$, for  $0\le j\le k_n$, if $n\ge N_1$.

\begin{lema}
\label{lem:4}
For all  $\epsilon>0$ there exists  $n\in \Z^+$ such that the
following holds:
\begin{enumerate}[i.]
\item $\dist(q_n,q^u) < \epsilon$.
\item If  $I$ is an  $E$-interval with  $q_n$ in the interior of  $I$
  and the length of one of each connected components of 
 $I\setminus \{q_n\}$ is greater than  $\epsilon$ then there exists
 $J\subset I$, $E$-interval, such that
$q_n\in\partial J$, $\ell(f^j(J)) \le \epsilon$, for  $0\le j \le
k_n$, and the extreme point of  $f^{k_n}(J)$ which is not
$f^{k_n}(q_n)$ belongs to the connected component of
$B_{\epsilon}(q^s)\cap W^s(p^s)$ that contains  $q^s$.
\end{enumerate}
\end{lema}

\begin{proof}
Let  $\psi:[-1,1]^2\to M$ be  $\Cr$ coordinates such that  $\psi(0,0)=q^s$;
$\eta\in C^E_\rho(x)$ (resp. $\eta\in C^F_\rho(x)$) is written as $(u,v)$
in such coordinates, with  $\|v\| < \|u\|$ (resp. $\|u\| < \|v\|$);
$B=\psi([-1,1]^2)\subset V_2$, and  $\Gamma=\psi([-1,1]\times\{0\})$
is the connected component of  $W^s(p^s)\cap B$ which contains  $q^s$.
 There exists $N_2\ge N_1$ such that  $f^{k_n}(q_n)\in
 \psi([-1/2,1/2]^2)$,
for all 
$n\ge N_2$. If  $I$ is an  $E$-interval with   $q_n\in\partial I$ we
shall say that it {\em point to}  $\Gamma$, if either there is   $J\subset I$, an
$E$-interval with   $q_n\in
J$ and  $f^{k_n}(J)\subset B$,  such that  the second
coordinate of every point in 
$f^{k_n}(J)$ is less than or equal to, in absolute value, the second
coordinate of  $f^{k_n}(q_n)$,
or $f^{k_n}(J)\cap \Gamma\neq \emptyset$.
Let  $\zeta_n = \sqrt{2} |y_n|$,
for  $n\ge N_2$, where $\psi(x_n,y_n)= f^{k_n}(q_n)$.
We observe that 
$\lim_{n\to\infty} \zeta_n = 0$, and if  $J$ is an  $E$-interval which
points to  $\Gamma$, with  $q_n\in\partial J$, and  $\ell(f^{k_n}(J))\ge \zeta_n$
then  $f^{k_n}(J)$ meets  $\Gamma$.

Let  $0<\epsilon<\min(\delta_0,\delta_1)$ be fixed and  $N_3\ge N_2$
such that 
$q_n\in B_\epsilon(q^u)$, for all  $n\ge N_3$, where  $\delta_0$ is
given by Proposition~\ref{pro:ps}. For each  $n\ge N_3$, let
$\tilde{I}_n$ be an  
$E$-interval that satisfies the hypothesis of the lemma,   i.e.
$\tilde{I}_n\setminus \{q_n\}$ has two connected components with its
length greater than  $\epsilon$. Let  $I_n$ be the unique 
$E$-interval contained  in  $\tilde{I}_n$ which point to  $\Gamma$ and
has lengths  $\epsilon$. We observe that by definition   $q_n\in\partial I_n$.
Let  $I_n^0, \ldots , I_n^{k_n}$ be the  $E$-intervals which satisfy:
$I_n^0\subset I_n$, $q_n\in \partial I_n^0$, $f^j(I_n^0) = I_n^j$ and
$\epsilon = \max\{\ell(I_n^j): 0\le j\le k_n\}$. We note that $\{I_n^j\}$
is well defined since $\ell(I_n)=\epsilon$.

If we prove that there is $n\ge N_3$ such that  $I_n^{k_n}$ intersects
 $\Gamma$
then we will have proved the lemma. Let us prove the existence of such
 an  $n$ by contradiction. So we assume that $I_n^{k_n}\cap \Gamma
= \emptyset$, for all  $n\ge N_3$, and get an absurd.

Let  $\nu_n = \max\{0\le j\le k_n: \ell(I_n^j)=\epsilon\}$. Note that 
$k_n - \nu_n\to \infty$, if  $n\to\infty$, since  $\zeta_n\to 0$.
Let $J_n =
I_n^{k_n}$, by compactness there exists a subsequence  $J_{n_j}$ which
converges to an 
 $E$-interval $J$ in such a way that  $x=\lim_{j\to\infty} f^{\nu_{n_j}}(q_{n_j})$
is well defined. It is clear that  $\ell(J)=\epsilon$ and  $x\in
\partial J$. Without loss of generality we assume that 
$J_n\to J$ and $f^{\nu_n}(q_n)\to x$.

 $J$ is an  $\epsilon$-$E$-interval (because of  $k_n-\nu_n\to\infty$) and
$x\in \Delta$, see \eqref{lemit:5} of  lemma~\ref{lem:conlem}. 
Proposition~\ref{pro:ps} implies that either $\omega(J)$ is a finite
collection
of invariant curves normally hyperbolic or  $\omega(J)$ is formed by
periodic points. Let us  see that each alternatives lead us to a
contradiction. We observe that 
   $\omega(J)\subset \La[0]$.

Consider the case when  $\omega(J)$ is a finite collection of invariant
 curves normally hyperbolic. Since the inner dynamics of theses curves
 are irrational and  $F$ restricted  to theses curves is its tangent
 space, the continuity of the dominated splitting implies that 
 $E$ is a normal fiber which is contractive. So 
$\omega(J)$ is an attractor, i.e. $W^s(\omega(J))$ is a neighborhood of
$\omega(J)$. Since  $x\in W^s(\omega(J))$, $q_n\in W^s(\omega(J))$ for $n$
big enough, but this is impossible, because  $\nu_n < k_n$, $f^{k_n}(q_n)\to q^s\in
W^s(p)$ and  $p\notin \omega(J)$.

Now consider the case when $\omega(J)$ is formed just by periodic
 points. Let  $y\in \omega(x)$, $\sigma$ and  $\tau$ be the
 eigenvalues associated to  $y$, with  $|\sigma| < |\tau|$, $W^-(y)$
 and  $W^+(y)$ be the invariant manifolds associated to  $\sigma$ and
 $\tau$ respectively. We observe that 
$|\sigma|\le 1$. With no loss of generality we suppose that 
 $y$ is a fixed point and $\sigma,\tau\in\R^+$. We have to consider
 several subcases. 

{\em $x$ belongs to the interior of  $W^s(y)$:} In this subcase we get
an absurd in an analogous way as when we consider $\omega(J)$ as a
finite collection of invariant curves normally hyperbolic. 

{\em $x$ does  not belong to the interior of  $W^s(y)$ and  $x= y$:} The proof
here is similar to the proof of  the next subcase. 

{\em $x$ does  not belong to the  interior of  $W^s(y)$ and $x\neq y$:}
Then  $\tau\ge
1$ and  $x\in W^-(y)$.  We get that  $W^+(y)$ has at most one stable
separatrix. Denote by  $\Sigma_x$ the separatrix of  $W^-(y)$ which
contains  $x$. We are going to obtain a contradiction.

%%%%%%%%%%%%%%%%%
\begin{figure}
\begin{center}
\includegraphics*{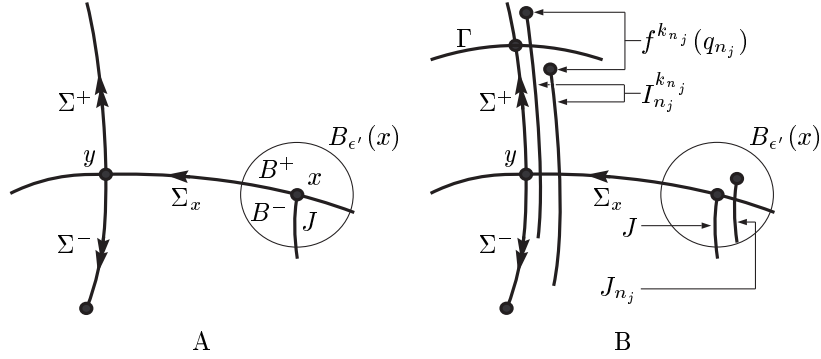}
\caption{Proof of Lemma  \ref{lem:4}}
\label{fig3}
\end{center}
\end{figure}
%%%%%%%%%%%%%%%%%

 Assume that 
$\tau>1$. Let  $\epsilon'>0$ be small and  $\Sigma_x^0$ be the
connected component of  $B_{\epsilon'}(x)\cap \Sigma_x$ which contains
 $x$.
$B_{\epsilon'}(x)\setminus \Sigma_x^0$ has two connected components
, we denote them by $B^-$ and  $B^+$. 
$J$ intercepts just one of them, let us  say 
 $B^-$. The  $\lambda$-lemma (\cite{palis}) implies that  $\{f^k(J)\}$
 accumulates in one of the separatrices of  $W^+(y)$, see
  Figure~\ref{fig3}-A.
Let $\Sigma^-$ be this separatrix  and $\Sigma^+$ be the other one of
$W^+(y)$. Since $\{f^k(J)\}$ accumulates on  $\Sigma^-$,
$\ell(\Sigma^-)\le \epsilon$.
\begin{afir}
If there are infinitely many  $n$ for which 
 $f^{\nu_n}(q_n)\in B^+$ then 
$\ell(\Sigma^+)\le \epsilon$.
\end{afir}

Let us prove the claim by contradiction. We have that $q^s\notin
\Sigma_x$. Let  $\{n_j\}$ be an increasing  sequence of positive
integers such that $f^{\nu_{n_j}}(q_{n_j})\in B^+$,  for all $j$, and
$n_1\ge N_3$. Since 
$\ell(\Sigma^+)>\epsilon$, there is  $r_j\in\Z^+$ such that
$\ell(f^{r_j}(J_{n_j}))>\epsilon$. From here it follows that 
 $\nu_{n_j}+r_j\ge
k_{n_j}$ and $q^s\in \Sigma^+$. Given  $j$ big enough we have two
alternatives. Either 
$I_{n_j}$ does not point to  $\Gamma$ or $I_{n_j}^{k_{n_j}}$
intersects
 $\Gamma$,
see Figure~\ref{fig3}-B. For both of the alternatives 
this is an absurd. This ends the proof of the Claim 1.

Suppose now that there is an increasing sequence of positive integers 
 $\{n_j\}$
such that  $f^{\nu_{n_j}}(q_{n_j})\in B^+$, for all  $j\ge 1$, then 
there exist  $z\in \Sigma^+$ and  $\{r_j\}$, with 
 $\nu_{n_j}<r_j<k_{n_j}$, such that 
$f^{r_j}(q_{n_j})\to z$. From here it follows that  $z\in\Delta$.
 Lemma~\ref{lem:3} give us the existence of a periodic point  $z'$ and
a small
neighborhood  $V$ of  $z$ which is contained in the stable manifold 
of the orbit of  $z'$, which is an absurd. We remark that the same 
conclusion is obtained if there is 
$\{n_j\}$ as above such that  $f^{\nu_{n_j}}(q_{n_j})\in B^-$. So we
have that 
$f^{\nu_n}(q_n)\in \Sigma_x$, for  $n$ big enough, again an absurd,
 because $\Sigma_x \subset W^s(y)$. With this we get the third subcase when  $\tau>1$.

If  $\tau=1$ then again there are several alternatives.
 $W^+(y)$ can have no separatrices, two or one;  moreover 
if it has one separatrix this can be an stable or an unstable
 separatrix. For each of this cases we get an absurd in the same way
 as above. This finishes the proof of the lemma.
\end{proof}

An immediate consequence of the lemma is the following result.

\begin{coro}
\label{cor:1}
If $\{\tilde{I}_n\}$ is a sequence of  $E$-intervals such that 
 $q_n\in
\tilde{I}_n$ and each of the connected components 
of  $\tilde{I}_n\setminus
\{q_n\}$ has length  $\delta_2=\min\{\delta_0,\delta_1\}/2$ then 
there are sequences $\{n_j\}\subset \Z^+$ and $\{J_j\}$ such that

\begin{enumerate}[i.]
\item $n_j\to\infty$
\item $J_j$ is an  $E$-interval with  $q_{n_j}\in \partial J_j$ and
$\ell(f^k(J_j))\le 1/j$, for $0\le k\le k_{n_j}$
\item $\partial f^{k_{n_j}}(J_j) = \{f^{k_{n_j}} (q_{n_j}), x_j\}$,
  with
$x_j\in \Gamma$, where $\Gamma$ is defined as in the lemma's proof. 
\end{enumerate}
\end{coro}

 Lemma~\ref{lem:4} has a version for  $F$-intervals, in such a  version
$f^{k_n}(q_n)$ is changed by  $q_n$ and  $f^{-1}$ by $f$.
We remark that the properties of $\{q_n\}$ used in the proof are:
(1) $\bigcup_{j=0}^{k_n}B_{\epsilon}(f^j(q_n))\subset \overline{V_1}$,
for
$\epsilon$ small and  $n$ big, and  (2) $\lim_{n\to\infty} q_n = q^u$
and 
$\lim_{n\to\infty} f^{k_n}(q_n) =q^s$. On the other hand,
$\{x_j\}$, given by  Corollary~\ref{cor:1}, satisfies: (1)
$\bigcup_{k=0}^{k_{n_j}} B_{\epsilon}(f^{-k}(x_j))\subset \overline{V_1}$,
if $\epsilon$ is small and $j$ big, and (2) $\lim_{j\to\infty} x_j = q^s$
and $\lim_{j\to\infty} f^{-k_{n_j}}(x_j) = q^u$. So a version for 
$F$-intervals is true changing the sequence  $\{q_n\}$ by the sequence
 $\{x_j\}$. In particular, if we take as  $F$-interval to $\Gamma$
we get that $W^s(p^s)$ intercepts transversally to $W^u(p^u)$ near of
$q^u$, so  $W^s(\Lambda)\cap W^u(\Lambda)\setminus \Lambda\neq
\emptyset$, which is a  contradiction, since we are assuming that
$W^s(\Lambda)\cap W^u(\Lambda)\setminus \Lambda = \emptyset$. 
So we conclude the proof of the theorem.\qed

%%%%%%%%%%%%%%%%%%%%%%%%%%%%%%%%%%%%%%%%%%%%%%%%%%%%%%%%%%%

\end{document}